\documentclass[a4paper]{amsart}
\usepackage[latin1]{inputenc}
\usepackage{ae,aecompl,amsbsy,amssymb,amsmath,amsthm,
eurosym,amsfonts,epsfig,graphicx,graphics,verbatim,enumerate,esint,color,MnSymbol}

\newcommand{\norm}[1]{\lVert #1 \rVert}

\newcommand{\br}{\mathbb{R}}

\newcommand{\cd}{\mathcal{D}}

\newcommand{\cs}{\mathcal{S}}

\newcommand{\dmu}{\mathrm{d}\mu}
\newcommand{\dnu}{\mathrm{d}\nu}
\newcommand{\dt}{\mathrm{d}t}

\usepackage{enumitem}

\theoremstyle{plain}% default
\newtheorem{theorem}{Theorem}[section]

\newtheorem{proposition}[theorem]{Proposition}

\theoremstyle{definition}
\newtheorem{definition}[theorem]{Definition}

\theoremstyle{remark}
\newtheorem*{remark}{Remark}

\numberwithin{equation}{section}

\begin{document}

%\subjclass[2010]{}
%    
\title{Equivalence of sparse and Carleson coefficients for general sets}
\author{Timo S. H\"anninen}
\address{Department of Mathematics and Statistics, University of Helsinki, P.O. Box 68, FI-00014 HELSINKI, FINLAND}
\email{timo.s.hanninen@helsinki.fi}
\begin{abstract}

We remark that sparse and Carleson coefficients are equivalent for every countable collection of Borel sets and hence, in particular, for dyadic rectangles, the case relevant to the theory of bi-parameter singular integrals. 

The key observation is that a dual refomulation by I. E. Verbitsky for Carleson coefficients over dyadic cubes  holds also for Carleson coefficients over general sets. We give a simple proof for this reformulation.
\end{abstract}
\date{\today}
\maketitle
%\tableofcontents
\section{}
The usual definitions of Carleson and sparse coefficients generalize word by word from the collection of dyadic cubes to an arbitrary countable collection of Borel sets:
\begin{definition}[Carleson coefficients in the generality of a collection of Borel sets]Let $\mu$ be a locally finite Borel measure on $\br^d$. Let $\cs$ be a countable collection of Borel sets. A family $\{\lambda_S\}_{S\in\cs}$ of non-negative reals is {\it Carleson} (with the constant $C\geq 1$) if we have
\begin{equation}\label{condition_equivalent}
\sum_{S\in\cs:S\subseteq \Omega} \lambda_S \leq C \mu(\Omega)
\end{equation}
for every union $\Omega$ of sets of the collection $\cs$.
\end{definition}
\begin{remark}
(a) It is equivalent to state the Carleson condition as follows: We have
\begin{equation}\label{condition_equivalent2}
\sum_{S\in \cs' } \lambda_S \leq C \mu(\bigcup_{S\in\cs'} S).
\end{equation}
for every subcollection $\cs'\subseteq \cs$. 

(b) For the collection $\cd$ of dyadic cubes, the equivalence between the seemingly stronger condition \eqref{condition_equivalent2} and the usual definition can be seen, for example, by first decomposing the collection $\cd'\subseteq \cd$ into the subcollections of cubes with a common ancestor and further splitting these subcollections into the subsubcollection of cubes with finitely many ancestors and into the subsubcollection of cubes with infinitely many ancestors and then applying the usual condition.

 (c) In the case of the Lebesgue measure and the collection of dyadic rectangles, it is equivalent to require the Carleson condition \eqref{condition_equivalent} to hold for every open set $\Omega$, in place of every union of dyadic rectangles, via approximating dyadic rectangles by open rectangles.

\end{remark}

\begin{definition}[Sparse coefficients in the generality of a collection of Borel sets]Let $\mu$ be a locally finite Borel measure on $\br^d$. Let $\cs$ be a countable collection of Borel sets.
 A family $\{\lambda_S\}_{S\in\cs}$ of non-negative reals is {\it sparse} (with the constant $C\geq 1$) if for each $S\in\cs$ there exists a subset $E_S\subseteq S$ such that $\lambda_S\leq C \mu(E_S)$ and such that the sets $\{E_S\}_{S\in\cs}$ are pairwise disjoint.
\end{definition}

The sparse coefficients are Carleson coefficients simply because
$$
\sum_{S:S\subseteq \Omega} \lambda_S\leq C \sum_{S:S\subseteq \Omega} \mu(E_S) \leq C \mu(\Omega)
$$
but the converse is more complicated, as well-known.
\begin{remark}Regarding the converse, note that Carleson coefficients may fail to be sparse coefficients by a simple obstruction: a point mass can not be divided. This obstruction is illustrated by the following example: Let $\delta_x$ be a Dirac measure at a point $x$. Then, for any nonzero coefficients $\lambda_{S_1}$ and $\lambda_{S_2}$ associated with any sets $S_1$ and $S_2$ such that $S_2\cap S_1 \ni x$, the coefficients  $\lambda_{S_1}$ and $\lambda_{S_2}$ are Carleson but unsparse. Therefore, the assumption that the measure $\mu$ has no point masses is needed in general for the converse.
\end{remark}

In the case of the collection of dyadic cubes, the converse was proven by I. E. Verbitsky \cite[Corollary 2]{verbitsky1996} by combining the following two steps:
\begin{itemize} \item First step: Dual reformulation of the Carleson condition as a certain estimate, by I. E. Verbitsky \cite[Theorem 4]{verbitsky1996}. This can be viewed as duality in discrete Littlewood--Paley spaces.
\item Second step: Characterizing such general estimates in terms of the existence of pairwise disjoint sets, by L. E. Dor \cite[Proposition 2.2]{dor1975}. This is based on functional and convex analysis.
\end{itemize}

In the case of the collection of dyadic cubes, an alternative proof (constructive hands-on proof avoiding functional analysis) was given by A. K. Lerner and F. Nazarov \cite[Lemma 6.3]{lerner2015} for the most important particular type of coefficients,  and this proof was generalized for the general type of coefficients by Cascante and Ortega \cite[Theorem 4.3]{cascante2017}.

The case of the collection of dyadic rectangles is relevant to the theory of bi-parameter singular integrals. Whether the converse holds in this case is mentioned as an open problem by A. Barron and J. Pipher \cite{barron2017} in their recent preprint. 

In this note, we observe that, thanks to the Dor--Verbitsky approach, the converse holds for the collection of dyadic rectangles and, more generally, for every countable collection of Borel sets:

\begin{theorem}[Carleson coefficients are sparse coefficients in the generality of a collection of Borel sets]Let $\mu$ be a locally finite Borel measure on $\br^d$. Assume that $\mu$ has no point masses. Let $\cs$ be a countable collection of Borel sets. Then, a family $\{\lambda_S\}_{S\in\cs}$ of non-negative reals is Carleson if and only if it is sparse, and moreover, the constants in both the conditions are the same.

\end{theorem}
In order to run the Dor--Verbitsky proof in this generality, we need only to prove the dual reformulation of the Carleson condition in this generality; an (elementary) proof for it is the contribution of this note. The proof is essentially a standard proof for the dyadic Carleson embedding theorem. 

In what follows, we explain the Dor--Verbitsky approach and prove the needed dual reformulation in this generality.
\subsection*{First step: Dual reformulation}

I.E. Verbitsky \cite[Theorem 4]{verbitsky1996} proves that the Carleson condition in the case of dyadic cubes can be rephrased as the following dual condition: The coefficients $\{\lambda_Q\}_{Q\in\cd}$ are Carleson if and only if we have the estimate
\begin{equation}\label{estimate_dual2}
\sum_{Q\in\cd} \lambda_Q a_Q \leq C \int \sup_{Q\in\cd} a_Q 1_Q \dmu
\end{equation}
for every family $a:=\{a_Q\}_{Q\in\cd}$ of nonnegative reals.

\begin{remark}The discrete Littlewood--Paley spaces $f^{p,q}(\mu)$ with $p,q\in(0,\infty]$ were essentially introduced by M. Prazier  and  B. Jawerth \cite{frazier1990}. For the exponents $p,q\in[1,\infty]$, the dual spaces $(f^{p,q}(\mu))^*=f^{p',q'}(\mu)$, where $p',q'$ denote the H\"older conjugate exponents of $p,q$, were computed by I.E. Verbitsky  \cite[Theorem 4]{verbitsky1996}. Therefore, in particular, the dual norm formula
\begin{equation}\label{estimate_littlewoodpaley}
\norm{\lambda}_{f^{\infty,1}(\mu)}=\sup\{\sum_{Q\in\cd} \lambda_Q a_Q \mu(Q) : a\in f^{1,\infty}(\mu)  \text{ with } \norm{a}_{f^{1,\infty}(\mu)}\leq 1\}
\end{equation}
 holds for the norms $$\norm{\lambda}_{f^{\infty,1}(\mu)}:=\sup_{Q\in\cd} \frac{1}{\mu(Q)} \sum_{\substack{R\in \cd: R\subseteq Q}} \lambda_R \mu(R) \quad \text{ and } \quad \norm{a}_{f^{1,\infty}(\mu)}:=\int \sup_{Q\in \cd} a_Q 1_Q \dmu. $$  
The dual norm formula \eqref{estimate_littlewoodpaley} for discrete Littlewood--Paley spaces states precisely the equivalence between the Carleson condition and the estimate \eqref{estimate_dual2}, and hence it is fitting to call the estimate \lq{}a dual reformulation\rq{} of the Carleson condition.
\end{remark}

The statement \eqref{estimate_dual2}  of the dual reformulation generalizes word by word to an arbitrary countable collection of Borels sets:
\begin{proposition}[Dual formulation in the generality of a collection of Borel sets]Let $\mu$ be a locally finite Borel measure on $\br^d$. Let $\cs$ be a countable collection of Borel sets. Let $\{\lambda_S\}_{S\in\cs}$ be a family of non-negative reals. Then, the following assertions are equivalent:
\begin{enumerate}[label=(\roman*)] 
\item \label{item_carleson} The family $\lambda$ is Carleson, that is $$
\sum_{S\in \cs'} \lambda_S \leq C \mu(\bigcup_{S\in\cs'} S)
$$
for every subcollection $\cs'$ of the collection $\cs$.
\item \label{item_estimate} We have the estimate
\begin{equation}\label{estimate_dual}
\sum_{S\in\cs} \lambda_S a_S \leq C \int \sup_{S\in\cs} a_S 1_S \dmu
\end{equation}
for every family $a:=\{a_S\}_{S\in\cs}$ of nonnegative reals.

\end{enumerate}
\end{proposition}
Next, we give an elementary proof for the statement. The proof is essentially a standard proof for the dyadic Carleson embedding theorem. 
\begin{proof}

First, we prove that the estimate \ref{item_estimate} implies the Carleson condition  \ref{item_carleson}. Let $\cs'$ be a subcollection of the collection $\cs$. We set $a_S:=1$ if $S\in \cs'$ and $a_S:=0$ otherwise, so that
\begin{equation}\label{temp_1}
\sum_{\substack{S\in\cs'}} \lambda_S=\sum_{S\in\cs} \lambda_S a_S. 
\end{equation}
By the assumed estimate \ref{item_estimate}, we have
\begin{equation}\label{temp_8}
\sum_{S\in\cs} \lambda_S a_S \leq C \int \sup_{S\in\cs} a_S 1_S \dmu.
\end{equation}
Observe that, by the choice of the coefficients $\{a_S\}_{S\in\cs}$, we have $\sup_{S\in\cs} a_S 1_S= 1_{ \bigcup_{S\in\cs'} S}$ and hence
\begin{equation}\label{temp_2}
\int \sup_{S\in\cs} a_S 1_S \dmu=\mu(\bigcup_{S\in\cs'} S).
\end{equation}
Combining the estimates \eqref{temp_1}, \eqref{temp_8}, and \eqref{temp_2} yields the claimed Carleson condition \ref{item_carleson}.

Next, we prove that the Carleson condition \ref{item_carleson} implies the estimate \ref{item_estimate}. Recall the distribution formula: For a non-negative measurable function on a measure space $(X,\nu)$, we have $\int_X f \dnu=\int_0^\infty \mu(f>t)\dt$. By applying this formula, we have
\begin{equation}\label{temp_3}
\sum_{S\in\cs} \lambda_S a_S=\int_{0}^ \infty \sum_{S\in \cs: a_S > t} \lambda_S \dt.
\end{equation}
Observe that $\{\sup_{S\in\cs} a_S 1_S>t\}=\bigcup_{\substack{S\in \cs: a_S > t}} S $. Therefore, by the assumed Carleson condition,
\begin{equation}\label{temp_4}
\sum_{S\in \cs: a_S > t} \lambda_S\leq  C \mu(\bigcup_{\substack{S\in \cs: a_S > t}} S)= C \mu(\{\sup_{S\in\cs} a_S 1_S>t\}).
\end{equation}
By applying the distribution formula again, we have
\begin{equation}\label{temp_5}
\int_{0}^ \infty \mu(\{\sup_{S\in\cs} a_S 1_S>t\})  \dt =\int \sup_{S\in\cs} a_S 1_S \dmu. 
\end{equation}
Combining the estimates \eqref{temp_3}, \eqref{temp_4}, and \eqref{temp_5} yields the claimed estimate \ref{item_estimate}.
\end{proof}

\subsection*{Second step: Characterization of the dual estimate via the existence of pairwise disjoint sets}
Estimates of the same form as the dual reformulation has been characterized by L. E. Dor \cite[Proposition 2.2]{dor1975} in terms of the existence of pairwise disjoint sets. Thanks to the (functional-analytic) generality of the characterization, the application of it to the case of Borel sets is word by word as in I. E. Verbitsky's \cite[Corollary 2]{verbitsky1996} application of it to the case of dyadic cubes. For the reader's convenience, we repeat the application here.

\begin{proposition}[Proposition 2.2 in L. E. Dor's article \cite{dor1975}]Let $\mu$ be a locally finite Borel measure on $\br^d$. Assume that $\mu$ has no point masses. Let $(g_i)_{i=1}^\infty $ be a sequence of non-negative measurable functions in $L^1(\mu)$. Let $C\geq 1$. Then, the following assertions are equivalent:
\begin{enumerate}
\item We have $$
\sum_{i=1}^\infty a_i \leq C  \int \sup_{i=1,2,\ldots} a_i g_i \dmu 
$$
for all sequences $(a_i)_{i=1}^\infty $ of nonnegative reals.
\item There exist disjoint measurable sets $E_1,E_2,\ldots$ in $\br^d$ such that
$$
1\leq C \int_{E_i} g_i \dmu
$$
for every $i=1,2,\ldots$.
\end{enumerate}
\end{proposition}
\begin{remark}In fact, L. E. Dor states and proves this proposition \cite[Proposition 2.2]{dor1975} for the Lebesgue measure on the unit interval $[0,1]$. In finding the auxiliary function $h$ in his proof \cite[Proof of Lemma 2.3]{dor1975}, he is implicitly using the property that the Lebesgue measure is {\it non-atomic}, which means that every set having positive measure can be splitted into two disjoint sets each having positive measure. Inspection of his proof shows that it works for every locally finite Borel measure $\mu$ on $\br^d$ that is non-atomic or, equivalently, that has no point masses.
\end{remark}

Recall that the dual reformulation states (after the change of variable $\tilde{a}_S:=\lambda_S a_S$) that
$$
\sum_{S\in\cs} \tilde{a}_S \leq C \int \sup_{S\in\cs} \tilde{a}_S \frac{1_S}{\lambda_S} \dmu
$$
for every family $\tilde{a}:=\{\tilde{a}_S\}_{S\in\cs}$ of nonnegative reals. Applying L E. Dor's characterization to the functions $g_S:= \frac{1_S}{\lambda_S}$ produces pairwise disjoint sets $\tilde{E}_S$ such that
$$
1\leq C \frac{1}{\lambda_S} \mu(\tilde{E}_S\cap S),
$$
and hence the sets $E_S:=\tilde{E}_S\cap S$ are the desired sets: these sets $E_S$ are pairwise disjoint, $E_S\subseteq S$, and $\lambda_S \leq C \mu(E_S)$.

{\bf Acknowledgement.} The author is grateful to Professor Tuomas Hyt\"onen for bringing to his attention the open problem of the equivalence between sparse and Carleson coefficients for rectangles, mentioned in A. Barron and J. Pipher's preprint \cite{barron2017}.
\bibliography{bibliography}{}

\begin{thebibliography}{1}

\bibitem{barron2017}
Alexander Barron and Jill Pipher.
\newblock Sparse domination for bi-parameter operators using square functions.
  {P}reprint.
\newblock 2017.
\newblock arXiv:1709.05009 [math.CA].

\bibitem{cascante2017}
Carme Cascante and Joaquin~M. Ortega.
\newblock Two-weight norm inequalities for vector-valued operators.
\newblock {\em Mathematika}, 63(1):72--91, 2017.

\bibitem{dor1975}
Leonard~E. Dor.
\newblock On projections in {$L_{1}$}.
\newblock {\em Ann. of Math. (2)}, 102(3):463--474, 1975.

\bibitem{frazier1990}
Michael Frazier and Bj\"orn Jawerth.
\newblock A discrete transform and decompositions of distribution spaces.
\newblock {\em J. Funct. Anal.}, 93(1):34--170, 1990.

\bibitem{lerner2015}
Andrei~K. Lerner and Fedor Nazarov.
\newblock Intuitive dyadic calculus: the basics. {B}ook.
\newblock 2015.
\newblock arXiv:1508.05639 [math.CA].

\bibitem{verbitsky1996}
Igor~E. Verbitsky.
\newblock Imbedding and multiplier theorems for discrete {L}ittlewood-{P}aley
  spaces.
\newblock {\em Pacific J. Math.}, 176(2):529--556, 1996.

\end{thebibliography}
\bibliographystyle{plain}
\end{document}